\documentclass[11pt,reqno]{amsart}
\usepackage{amsfonts,amsmath,amssymb}
\def\beq#1#2\eeq{%
        \begin{equation}%
        \label{#1}%
            #2%
        \end{equation}%
    }

\usepackage{color}
\usepackage{graphicx}
\usepackage{epstopdf}

\theoremstyle{plain}
\newtheorem{theorem}{Theorem}

\theoremstyle{remark}

\theoremstyle{definition}

\def\a{\alpha}

\title[Euler-Dirac)]{New integrable two-centre problem on sphere in Dirac magnetic field}

\author{A.P. Veselov}
\address{Department of Mathematical Sciences,
Loughborough University, Loughborough LE11 3TU, UK; Moscow State University and Steklov Mathematical Institute, Moscow, Russia}
\email{A.P.Veselov@lboro.ac.uk}

\author{Y. Ye}
\address{Department of Mathematical Sciences,
Loughborough University, Loughborough LE11 3TU, UK}
\email{Y.Ye@lboro.ac.uk}

\begin{document}

\maketitle

\begin{abstract}
We present a new family of integrable versions of the Euler two-centre problem on two-dimensional sphere in the presence of the Dirac magnetic monopole of arbitrary charge. The new systems have very special algebraic potential and additional integral quadratic in momenta, both in classical and quantum versions. \end{abstract}

%\tableofcontents

\section{Introduction}

The celebrated Euler two-centre problem \cite{Euler} was one of the first non-trivial mechanical systems integrated completely since the solution of famous Kepler problem by Newton.
In its two-dimensional version the Hamiltonian has the form
$$H=\frac{1}{2}(p_1^2+p_2^2)-\frac{\mu}{r_1}-\frac{\mu}{r_2},$$
where $$r_1=\sqrt{(q_1+c)^2+q_2^2}, \, r_2=\sqrt{(q_1-c)^2+q_2^2}$$ are the distances from the two centres fixed at the points $(\pm c,0).$

In the confocal coordinates $u_1=r_1+r_2, u_2=r_1-r_2$ the variables in the corresponding Hamilton-Jacobi equation can be separated, leading to the explicit solution of the system in quadratures (see Arnold \cite{Arnold}). More recent detailed analysis of this classical system can be found in Waalkens et al \cite{WDR}.

Its natural generalisation to the spaces of constant curvature was originally found in 1885 by Killing \cite{K} and rediscovered by Kozlov and Harin \cite{KH} as part of a general family of systems, separable in spherical elliptic coordinates. More of the history of this problem with the references can be found in Borisov and Mamaev \cite{BM}, who also discussed various integrable generalisations of this system.

The main observation (in the hyperbolic case due already to Bolyai and Lobachevsky and in spherical case to Serret (see \cite{BM}) is that for the system on the unit sphere $S^2$ the non-Euclidean analogue of the Newtonian potential $\mu/r$ is $\mu\cot \theta$, where $\theta$ is the spherical distance between the particle and fixed centre. 
The dynamics of corresponding natural 2-centre version with 
$$
U=-\mu \cot \theta_1-\mu\cot \theta_2,
$$
where $\theta_1$ and $\theta_2$ are the spherical distances from the fixed centres, was first studied by Killing \cite{K}, who separated the variables in the corresponding Hamilton-Jacobi equation using the Neumann elliptic coordinates on sphere.

Note that this potential $U$ has actually {\it four singularities} on the sphere, which can be interpreted as two antipodal pairs of centres with opposite charges $\pm \mu$ (see Fig. 1),  and thus this system should probably be considered as a 4-centre problem with Coulomb, rather than gravitational Newtonian interaction.

\begin{figure}[h]
  \includegraphics[width=54mm]{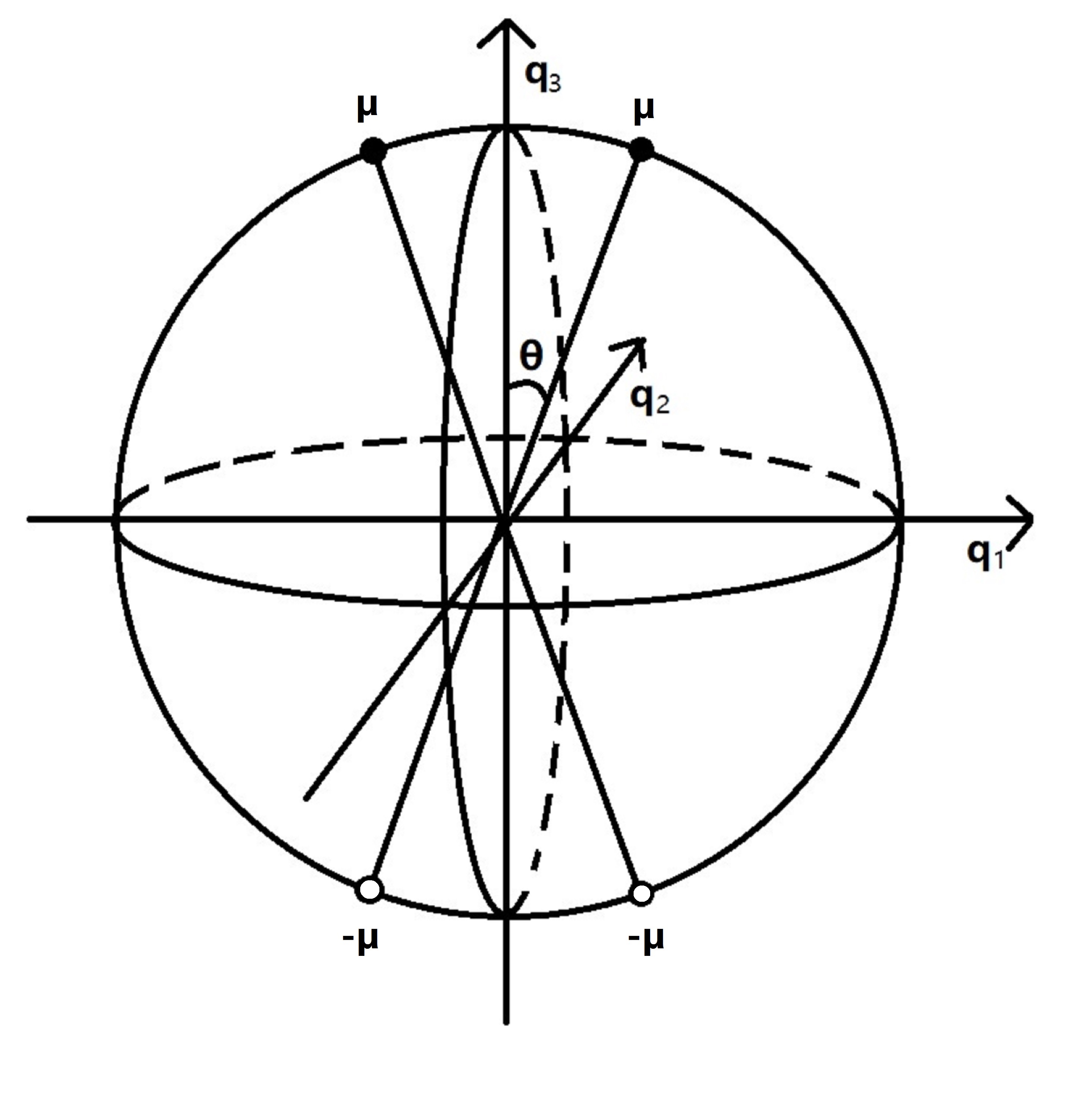} \quad   \includegraphics[width=50mm]{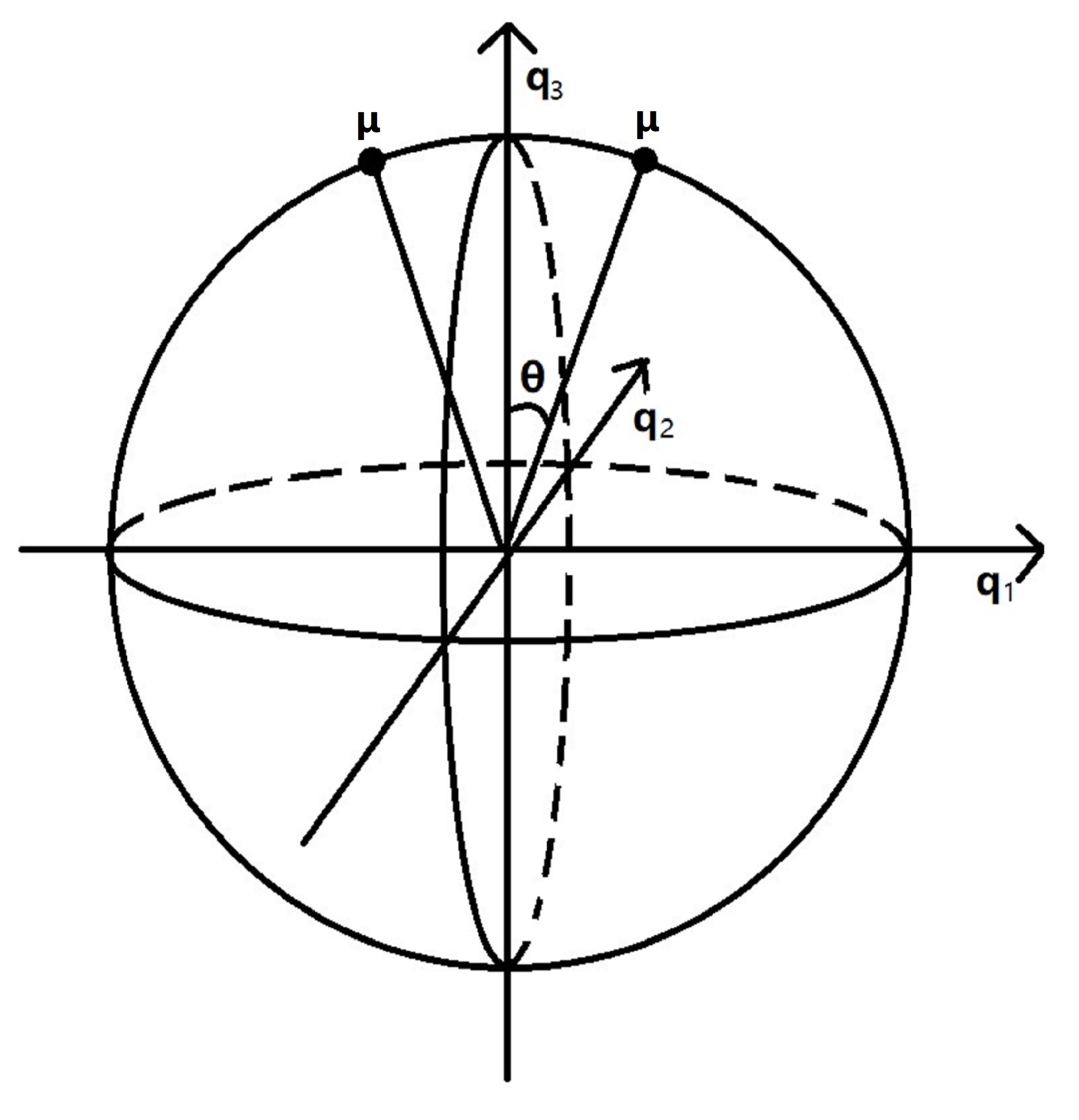}    
  \caption{Position of the fixed centres in the classical (left) and the new (right) systems}
  \end{figure}
  
To write down the explicit formula for the Hamiltonian and the additional integral it is convenient to use the canonical Lie-Poisson bracket on the dual space $e(3)^*$ of the Lie algebra of the group of motion of the Euclidean space. The corresponding variable $M_i, q_i, \, i=1,2,3$ have the Poisson brackets
\begin{equation}
\label{e3}
\left\{M_i, M_j\right\} = \epsilon_{ijk} M_k, \ \  \left\{M_i, q_j\right\} = \epsilon_{ijk} q_k, \ \  \left\{q_i, q_j \right\}=0.
\end{equation}
We have two Casimir functions 
$$
C_1=|q|^2,\quad C_2=(M,q).
$$
The symplectic leaf with $C_1=|q|^2=1, C_2=(M,q)=0$ is symplectically isomorphic to the cotangent bundle of the unit sphere $T^*S^2.$

In the coordinates $M,q$ the Hamiltonian of the spherical analogue of the Euler two-centre problem is
%\begin{equation}
%\label{Hold}
$$
H=\frac{1}{2}|M|^2-\mu \frac {\beta q_3-\alpha q_1}{\sqrt {q_2^2+(\alpha q_3+\beta q_1)^2}}-\mu \frac {\beta q_3+\alpha q_1}{\sqrt {q_2^2+(\alpha q_3-\beta q_1)^2}},
$$
%\end{equation}
where $\mu,\alpha,\beta$ are parameters such that $\alpha^2+\beta^2=1.$ 

We have 4 fixed centres at $(\pm \alpha, 0,\pm\beta)$, two of which at $(\pm \alpha, 0,\beta)$ for $\mu>0$ are attractive, while 
their antipodes $(\pm \alpha, 0,-\beta)$ are repulsive (see Fig. 1).

The explicit form of the additional integral at the special Casimir level $(M,q)=0$ was found by Mamaev \cite{M} (see also \cite{BM})\footnote{There is a couple of sign typos in the explicit form of the integral in these papers, which were kindly corrected for us by Ivan Mamaev.}:
%\begin{equation}
%\label{Fold}
$$
F=\alpha^2 M_1^2-\beta^2 M_3^2-2\alpha\beta\left(\mu \frac {\beta q_1-\alpha q_3}{\sqrt {q_2^2+(\beta q_1-\alpha q_3)^2}}+\mu \frac {\alpha q_1+\beta q_3}{\sqrt {q_2^2+(\alpha q_3+\beta q_1)^2}}\right).
$$
%\end{equation}
For the recent detailed analysis of the orbits in this system see the paper by Gonzalez Leon et al \cite{GGM}.

Note that on the symplectic leaves with $(M,q)\neq 0$ (corresponding to the additional non-zero Dirac magnetic field) the system is believed to be non-integrable for non-zero values of parameters. From \cite{VY} it follows that at least it has no additional integrals quadratic in momenta.

The aim of this paper is to introduce a new integrable family of systems on $e(3)^*$ with the Hamiltonian
$$
H=\frac{1}{2}|M|^2-\mu\frac{|q|}{\sqrt{R(q)}}, \,\,\,\,\, R(q)=Aq_2^2+Bq_1^2+(A+B)q_3^2-2\sqrt{AB}|q|q_3
$$
depending on real parameters $\mu, A, B$ satisfying $A>B>0$. 
These systems can be interpreted as the motion on the unit sphere with very particular algebraic electric potentials having two singularities of Coulomb/Newtonian type in the external field of Dirac magnetic monopole (see more in the next section). 

We show that the new systems are integrable in both classical and quantum case for all values of parameters and {\it for all values} of magnetic charge (or, equivalently, for all values of Casimir function $C_2=(M,q)$). The additional integral
$$
F=A M_1^2+B M_2^2+\frac{2 \sqrt{A B}}{|q|}(M,q) M_3-2\mu\sqrt{AB}\frac{q_3}{\sqrt{R(q)}}
$$
is quadratic in $M$ but with coefficients depending algebraically on $q$. 

The significance of the new systems is explained by the fact that this is the only integrable family of this type apart from the classical Clebsch system, which can be interpreted as the Neumann system on sphere with quadratic potential and additional Dirac magnetic field (see \cite{V,VY}).

 \section{New system}
 
Consider now the general symplectic leaves in $e(3)^*$, which are the coadjoint orbits of the Euclidean motion group $E(3)$ determined by
\begin{equation}
\label{symleaf}
(q,q)=R^2, \qquad (M,q)=\nu R.
\end{equation} 
Novikov and Schmelzer \cite{NS} introduced the variables
\begin{equation}
\label{NS}
 L_i=M_i-\frac{\nu}{R}q_i, \quad i=1,2,3
 \end{equation}
to identify the coadjoint orbits with $T^*S^2:$ 
$
(q,q)=R^2,\, (L,q)=0.
$

Assuming here for the convenience that the radius of sphere is 1, we have
in the new variables the Poisson brackets are
\begin{equation}
\left\{L_i, L_j\right\} = \epsilon_{ijk}\left(L_k -\nu q_k \right), \ \  \left\{L_i, q_j\right\} = \epsilon_{ijk} q_k, \ \  \left\{q_i, q_j \right\}=0 \label{eq:PBs}
\end{equation}
and the corresponding symplectic form becomes
\begin{equation}
\omega=\, \mathrm{d} P \wedge \, \mathrm{d} Q +  \nu  \, \mathrm{d} S, \label{eq:symp}
\end{equation}
where $\, \mathrm{d} P \wedge \, \mathrm{d} Q$ is the standard symplectic form on $T^* S^2$ and $\, \mathrm{d} S$ is the area form on $S^2$ (see \cite{NS}). 
As it was pointed out in \cite{NS} the second term corresponds to the magnetic field of the Dirac monopole:
$
\mathcal H= \nu \, \mathrm{d} S.
$

 Consider the following functions on $e(3)^*$:
\begin{equation}
\label{H}
H=\frac{1}{2}|M|^2-\mu\frac{|q|}{\sqrt{R(q)}},
\end{equation}
\begin{equation}
\label{F}
F=A M_1^2+B M_2^2+\frac{2 \sqrt{A B}}{|q|}(M,q) M_3-2\mu\sqrt{AB}\frac{q_3}{\sqrt{R(q)}},
\end{equation}
where
\begin{equation}
\label{R}
R(q)=Aq_2^2+Bq_1^2+(A+B)q_3^2-2\sqrt{AB}|q|q_3
\end{equation}
and $\mu, A, B$ are parameters satisfying $A>B>0$.

The graph of the potential $U=-\mu\frac{|q|}{\sqrt{R(q)}}$ after stereographic projection from the North pole is shown on Fig. 2.
We can interpret it as the Coulomb-like potential with two fixed charged centres (repulsive when $\mu<0,$ and attractive when $\mu>0$).

\begin{figure}[h]
  \includegraphics[width=80mm]{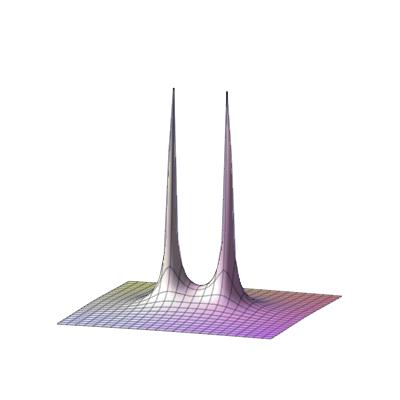}   % \includegraphics[width=60mm]{-potential2}   
  \caption{Graph of the repulsive version of the potential $U$ after stereographic projection}
  \end{figure}

Note that our system with positive $\mu$ admits also Newtonian gravitational interpretation in contrast to the classical Killing version.

\begin{theorem} 
The Poisson bracket of $H$ and $F$ is identically zero: $$\{F,H\}\equiv 0.$$ On the symplectic leaf with $|q|^2=1, (M,q)=\nu$ we have a new integrable system on unit sphere with the potential having two Coulomb-like singularities with charge $\mu/\sqrt{A-B}$ fixed at the points $(\pm \alpha,0, \beta)$ with
$$
\alpha=\sqrt{\frac{A-B}{A}}, \quad \beta=\sqrt{\frac{B}{A}}
$$
 in the external field of Dirac magnetic monopole with charge $\nu$.
 \end{theorem}

The proof can be done by direct calculation. The formulae (\ref{H}),(\ref{F}),(\ref{R}) was the final outcome of the lengthy local investigations started many years ago by the first author in collaboration with Ferapontov and Sayles, see \cite{VY}.

Let us check that the potential $U(q)=-\mu/\sqrt{R(q)}$ has the described singularities. For this re-write $R(q)$ as
$$
R(q)=(A-B)q_2^2+(\sqrt{A} q_3-\sqrt{B} |q|)^2.
$$
Indeed, we have
$$
(A-B)q_2^2+(\sqrt{A} q_3-\sqrt{B} |q|)^2=(A-B)q_2^2+A q_3^2-2\sqrt{AB}|q|q_3 +B (q_1^2+q_2^2+q_3^2)
$$
$$
=Aq_2^2+Bq_1^2+(A+B)q_3^2-2\sqrt{AB}|q|q_3.
$$
Since $A>B>0$ the equality $$(A-B)q_2^2+(\sqrt{A} q_3-\sqrt{B} |q|)^2=0$$ implies that both $q_2=0$ and $\sqrt{A} q_3-\sqrt{B}|q|=0,$ which gives two points
$$
q_1=\pm\sqrt{\frac{A-B}{A}}|q|, \, q_2=0,\, q_3=\sqrt{\frac{B}{A}}|q|.
$$
One can check that near the singularities we have
$$
R \approx 2(A-B) \rho^2,
$$
where $\rho<<1$ is the spherical distance from the singularity, which means that
$$
U \approx -\frac{\mu}{\sqrt{2(A-B)}\rho}
$$
is Coulomb-like as claimed.

Note that at the antipodal points $(\pm \alpha,0,-\beta)$ the potential is smooth since 
$$
R(\pm \alpha,0,-\beta)=-2\sqrt B \neq 0.
$$

The Euler equations on $e(3)^*$ are known to coincide with the classical Kirchhoff equations of the rigid body in the infinite ideal fluid \cite{NS}.
The only other known integrable case with the Hamiltonian
$H=\frac{1}{2}|M|^2+U(q)$ is the classical Clebsch system \cite{Clebsch} with quadratic potential
$$
U=\frac{1}{2}(a_1q_1^2+a_2q_2^2+a_3q_3^2).
$$
One can show that this these two systems exhaust all integrable cases in this class with quadratic in $M$ integral \cite{VY}.

In the second known classical integrable Steklov-Lyapunov case \cite{Lyapunov, Steklov} of the Kirchhoff equations the Hamiltonian 
$$
H=\frac{1}{2}|M|^2-(M,Aq)+\frac{1}{2}(C^2q,q) 
$$
with $\,\, A=diag\,(a_1, a_2, a_3), \,\, C=diag\,(a_2-a_3, a_3-a_1, a_1-a_2),$ has the terms linear in $M$, responsible for the additional non-Dirac magnetic field.

\section{New system in spherical elliptic coordinates}

Consider the unit sphere given in the Cartesian coordinates $q_1$, $q_2$, $q_3$ by the equation
\begin{equation}
  q_1^{\,2}+q_2^{\,2}+q_3^{\,2}=1, \nonumber
\end{equation}
and introduce the {\it spherical elliptical coordinates} \cite{Moser,N} as the roots $u_1$, $u_2$ of the quadratic equation
\begin{equation}
\label{elli}
  \Phi(u)=\frac{q_1^{\,2}}{A-u}+\frac{q_2^{\,2}}{B-u}+\frac{q_3^{\,2}}{C-u}=0, 
\end{equation}
\\[-4mm]
where $C=0$ and $A>B>0$ are the same as before.  
Expressing the function
\begin{equation}
  \Phi(u) =\frac{(u-u_1)(u-u_2)}{-u(A-u)(B-u)}
\end{equation}
 in terms of partial fractions we have the formula for the Cartesian coordinates in terms of the elliptic coordinates
\begin{equation}
\label{qs}
  q_1^{\,2}= \frac{(A-u_1)(A-u_2)}{A(A-B)}, \ \ \
  q_2^{\,2}= \frac{(B-u_1)(B-u_2)}{B(B-A)}, \ \ \
  q_3^{\,2}= \frac{u_1u_2}{AB}. 
\end{equation}
The metric on the sphere in the elliptic coordinates takes the form
\begin{equation}
\label{metrics}
  ds^2=  \frac{u_1-u_2}{f(u_1)} du_1^2+  \frac{u_2-u_1}{f(u_2)} du_2^2.
\end{equation}
where 
\begin{equation}
\label{f}
f(u)=-4u(u-A)(u-B).
\end{equation}
The graph of $f(u)$ and the positions of the variables $$0\leq u_1\leq  B\leq u_2\leq A$$ are shown on Fig. 3.

\begin{figure}[h]
  \includegraphics[width=65mm]{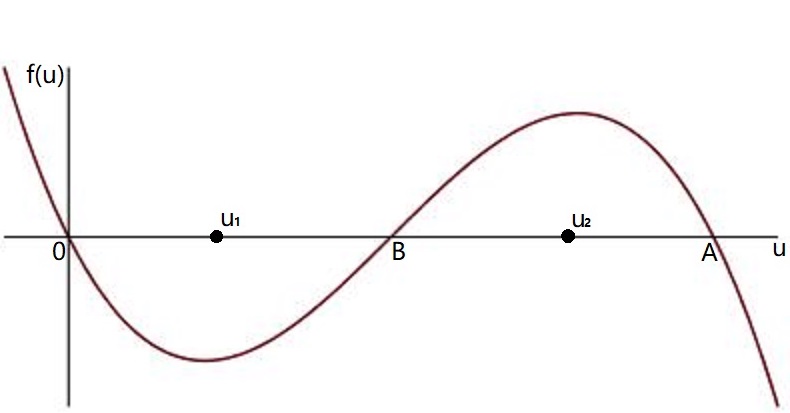}   
  \caption{Choice of roots $u_1,u_2$}
  \end{figure}
  
From (\ref{elli}) we have
$$
(u-u_1)(u-u_2)=u^2-(Bq_1^2+Aq_2^2+(A+B)q_3^2)u+ABq_3^2,
$$
implying
\begin{equation}
  u_1+u_2= Bq_1^2+Aq_2^2+(A+B)q_3^2, \quad u_1u_2=ABq_3^2.
\end{equation}
In particular,
$$
(\sqrt{u_1}-\sqrt{u_2})^2=u_1+u_2-2\sqrt{u_1u_2}=Bq_1^2+Aq_2^2+(A+B)q_3^2-2\sqrt{AB}q_3
$$
$$
=B(q_1^2+q_2^2+q_3^2)+ (A-B)q_2^2+Aq_3^2-2\sqrt{AB}q_3=(A-B)q_2^2+(\sqrt{A}q_3-\sqrt{B})^2,
$$
so 
$$
R(q)=(A-B)q_2^2+(\sqrt{A}q_3-\sqrt{B})^2=(\sqrt{u_1}-\sqrt{u_2})^2.
$$
Thus the (electric) potential of our new system is 
$$
U(q)=-\frac{\mu}{\sqrt{R(q)}}=-\frac{\mu}{\sqrt{u_2}-\sqrt{u_1}}.
$$
Note that we have chosen here $q_3=\sqrt{u_1u_2}/\sqrt{AB}$ to be positive, which leads to the singularities corresponding to
$
u_1=u_2=B
$ 
lying in the upper half-space.

To write down the kinetic term we need to introduce the potential of the Dirac magnetic monopole on sphere, which is known to be impossible to choose non-singular for non-zero magnetic charge $\nu$ by topological reasons (see e.g. Wu-Yang \cite{WY}). If we make two punctures at North and South poles, then we can use, for example, 
\begin{equation}
\label{A}
A=\nu q_3\frac{q_1dq_2-q_2dq_1}{q_1^2+q_2^2}.
\end{equation}
Indeed, one can check that on the sphere $|q|^2=1$ 
$$
dA=\nu(q_1dq_2\wedge dq_3+q_2dq_3\wedge dq_1+q_3dq_1\wedge dq_2)=\nu dS
$$
(where $dS$ is the area form on the unit sphere), which is the magnetic form of the Dirac monopole with charge $\nu.$

Let $$A=A_1(u)du_1+A_2(u)du_2$$ be any such 1-form (e.g. given by (\ref{A})) written in the elliptic coordinates, so that
$$
dA=\nu dS=B(u)du_1\wedge du_2,
$$
where
\begin{equation}
\label{B}
 B(u)=\nu\frac{u_2-u_1}{\sqrt{-f(u_1)f(u_2)}}
\end{equation}
is the density of the Dirac magnetic field in the elliptic coordinates.

Define the {\it magnetic momenta} by
\begin{equation}
\label{pi}
\tilde p_i=p_i-A_i(u), \, i=1,2.
\end{equation}
The corresponding Poisson brackets are
$$
\{\tilde p_1,\tilde p_2\}=B(u), \, \{\tilde p_1,u_1\}=\{\tilde p_2,u_2\}=1,
$$
with all other to be zero.

The Hamiltonian and the integral of the new system can be written now as
\begin{equation}
\label{H1}
H=\frac{1}{2}\left(\frac{f(u_1)}{u_1-u_2}\tilde p_1^2+\frac{f(u_2)}{u_2-u_1}\tilde p_2^2\right)-\frac{\mu}{\sqrt{u_2}-\sqrt{u_1}}
\end{equation}
\begin{equation}
\label{F1}
F=u_2\frac{f(u_1)}{u_1-u_2}\tilde p_1^2+u_1\frac{f(u_2)}{u_2-u_1}\tilde p_2^2+\phi_1(u) \tilde p_1+\phi_2(u) \tilde p_2+V(u),
\end{equation}
where
$$
\phi_1=-\nu\frac{\sqrt{-f(u_1)f(u_2)}}{\sqrt{u_1u_2}+u_2}, \, \phi_2=-\nu\frac{\sqrt{-f(u_1)f(u_2)}}{\sqrt{u_1u_2}+u_1}
$$
and
$$
V=-\frac{2\mu\sqrt{u_1u_2}}{\sqrt{u_2}-\sqrt{u_1}}-\nu^2 (\sqrt{u_1}-\sqrt{u_2})^2.
$$
Note that the electric potential can be written in St\"ackel form as
$$
U=-\frac{\mu}{\sqrt{u_2}-\sqrt{u_1}}=-\frac{\mu(\sqrt{u_2}+\sqrt{u_1})}{u_2-u_1},
$$
so, when the magnetic charge $\nu=0$, we can take $A=0,$ $\tilde p_i=p_i$ and the variables in the corresponding Hamilton-Jacobi equation 
$$
\frac{f(u_1)}{u_1-u_2}\left(\frac{\partial S}{\partial u_1}\right)^2+\frac{f(u_2)}{u_2-u_1}\left(\frac{\partial S}{\partial u_2}\right)^2-\frac{\mu(\sqrt{u_2}+\sqrt{u_1})}{u_2-u_1}=h
$$
can be separated (similarly to the classical two-centre problem, see Kozlov-Harin \cite{KH}).

What makes our case special is that the integrability holds for general $\nu$, although the separation of variables does not work, at least immediately (see Concluding remarks).

Another nice property of the new system is that it has a natural integrable quantum version.

\section{Quantum version}

Let us first recall the geometric quantisation of the Dirac magnetic monopole following Kemp and one of the authors \cite{KV}.

Let $$X_1=q_3 \partial_2 - q_2 \partial_3,\ X_2= q_1 \partial_3 - q_3 \partial_1,\ X_3= q_2 \partial_1 - q_1 \partial_2$$ be the vector fields generating rotations of $S^2$ given by $q_1^2+q_2^2+q_3^2=1$ and $\nabla_{X_j}$ be the corresponding covariant derivatives with respect to the Dirac $U(1)$-connection $iA$ with $A$ given, for example, by (\ref{A}).

Note that in the quantum case the charge $\nu$ of the Dirac magnetic monopole must be quantised, as it was pointed out already by Dirac \cite{Dirac}.
Geometrically this corresponds to the integrality of the Chern class of $U(1)$-bundle over sphere:
$$
\frac{1}{2\pi} \int_{S^2} \nu \mathrm{d} S=2\nu \in \mathbb Z.
$$

Then one can check \cite{KV} that
$$\hat \nabla_j:=i \nabla_{X_j}$$ and the operators $\hat q_j$ of multiplication by  $q_j$ satisfy the commutation relations
$$[\hat \nabla_k, \hat \nabla_l] = i \epsilon_{klm}(\hat \nabla_m - \nu \hat q_m)$$
and thus can be considered as quantization of Novikov-Schmelzer variables.
The quantum versions of the original variables
\begin{equation}
\hat M_j = \hat \nabla_j + \nu q_j,
\label{mam}
\end{equation}
satisfy the standard angular momentum relations 
\begin{equation*}
[\hat M_k, \hat M_m] = i\epsilon_{kmn}\hat M_n, \, [\hat M_k, \hat q_m] = i\epsilon_{kmn}\hat q_n,
\end{equation*} 
and coincide with Fierz's modification of the angular momentum 
in the presence of the Dirac magnetic monopole \cite{Fierz}. 

The quantum Hamiltonian of the Dirac monopole can be written 
 in terms of magnetic angular momentum $\hat M$ as
\begin{equation*}
\hat H=\frac{1}{2}(\hat M_1^2+\hat M_2^2+\hat M_3^2).
\end{equation*}
Since the operator $\hat H$ is a Casimir operator for $SO(3)$ it 
acts on every irreducible representation of $SO(3)$ as a scalar, which allows to use the representation theory for an explicit computation of the spectrum of Dirac magnetic monopole on the Hilbert space of functions $\mathcal F(S^2)$ (see the details in \cite{KV}).
 
The quantum Hamiltonian of our system can now be defined by the same formula 
\begin{equation}
\label{qH}
\hat H=\frac{1}{2}(\hat M_1^2+\hat M_2^2+\hat M_3^2)-\mu\frac{|q|}{\sqrt{R(q)}},
\end{equation}
where as before
$
R(q)=(A-B)q_2^2+(\sqrt{A} q_3-\sqrt{B} |q|)^2.
$
Here slightly abusing notations, we mean by $f(q)$ the operator of multiplication by $f(q)$ on the space of functions $\mathcal F(S^2)$.
Define similarly the quantum version of $F$ by 
\begin{equation}
\label{qF}
\hat F=A \hat M_1^2+B \hat M_2^2+ \frac{2 \sqrt{A B}}{|q|}(\hat M,q)\hat M_3-2\mu\sqrt{AB}\frac{q_3}{\sqrt{R(q)}}.
\end{equation}
Note that there is no ordering problem here since $\hat M_3$ commutes with the operator $(\hat M,q)|q|^{-1}.$

\begin{theorem} 
Operators $\hat H$ and $\hat F$ given by (\ref{qH}), (\ref{qF}) commute: 
$$
[\hat H, \hat F]=0,
$$ 
so $F$ is the second quantum integral of the system, ensuring its integrability.
\end{theorem}

The proof is again by straightforward check. When the magnetic charge $\nu=0$ the Hamiltonian and integral become
$$
\hat H=-\frac{1}{2}\Delta -\mu\frac{|q|}{\sqrt{R(q)}}, \quad \hat F=A \hat M_1^2+B \hat M_2^2-2\mu\sqrt{AB}\frac{q_3}{\sqrt{R(q)}},
$$
where $\Delta$ is the Laplace-Beltrami operator on the unit sphere.

\section{Hyperbolic version}

Replacing the Euclidean group of motion $E(3)$ by the group $E(2,1)$ of motion of the pseudo-Euclidean space $\mathbb R^{2,1}$ we come to the following natural hyperbolic version of our system.

The corresponding Lie algebra of $so(2,1)$ consists of $3\times 3$ matrices $X$ satisfying $X:XJ+JX^T=0,$ where 
$$J=\left( \begin{array}{ccc}
	-1&0&0\\0&-1&0\\0&0&1
\end{array}\right).$$
We have a natural basis $M_1,M_2,M_3,$ defined by 
$$X=\left( \begin{array}{ccc}
0&-a_3&a_2\\a_3&0&a_1\\a_2&a_1&0
\end{array}\right)=a_1M_1+a_2M_2+a_3M_3,$$
with the commutation relations between themselves and with natural generators of translations $q_1,q_2,q_3$:
\begin{equation}
\label{hP}
[M_1,M_2]=M_3,\,[M_2,M_3]=-M_1,\,[M_3,M_1]=-M_2,
\end{equation}
$$[M_1,q_2]=q_3=-[M_2,q_1],\, [M_1,q_3]=q_2=-[M_3,q_1],\, [M_2,q_3]=-q_1=-[M_3,q_2].$$

The Casimir functions are $$C_1=(q, Jq)=-q_1^2-q_2^2+q_3^2:=||q||^2,$$
defining the pseudo-Euclidean structure on $\mathbb R^{2,1}$, and
$$C_2=  \langle M,q\rangle:=(M, Jq)=-M_1q_1-M_2q_2+M_3q_3.$$
The relation $C_1=||q||^2=1$ now defines the two-sheeted hyperboloid, one sheet of which presenting a model of the hyperbolic plane (see Fig. 4).

\begin{figure}[h]
  \includegraphics[width=60mm]{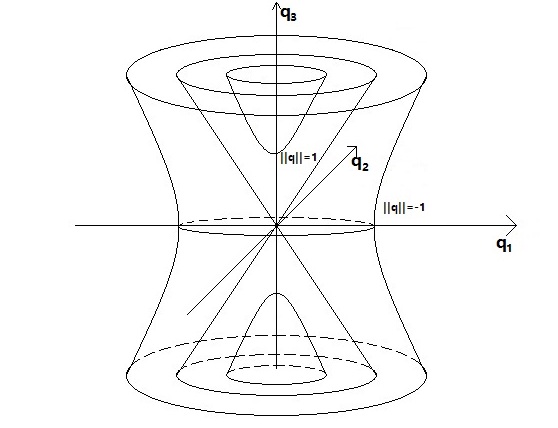}   
  \caption{Pseudo-spheres in $\mathbb R^{2,1}$}
  \end{figure}

The Hamiltonian of the natural hyperbolic analogue of the new system is
\begin{equation}
\label{hH}
H=\frac12(M_1^2+M_2^2-M_3^2)+\frac{\mu ||q||}{\sqrt{R(q)}},
\end{equation}
where
\begin{equation}
\label{hR}
R(q)=-A{{q_2}}^{2}-B{{ q_1}}^{2}+ \left( A+B \right) q_3^{2}-2\,\sqrt {AB}||q||q_3,
\end{equation}
which can be rewritten as
$$
R(q)=(B-A)q_2^2+(\sqrt{A}q_3-\sqrt{B}||q||)^2.
$$
Assuming now that $B>A>0,$ we see that the potential has two singularities when $q_2=0, \, \sqrt{A}q_3-\sqrt{B}||q||=0,$ or, if we assume that $||q||^2=1$, at two points
$
(\pm \frac{\sqrt{B-A}}{\sqrt{A}}, 0, \frac{\sqrt{B}}{\sqrt{A}}).
$

The corresponding additional integral has the form
\begin{equation}
\label{hF}
F=A{{ M_1}}^{2}+B{{ M_2}}^{2}-2\,{\frac {\sqrt {AB}}{||q||} \langle M,q \rangle M_3}+2\,{\frac {\mu\,\sqrt {AB}{ q_3}}{\sqrt{R(q)}}}.
\end{equation}

\begin{theorem} 
The functions $H$ and $F$ given by (\ref{hH})-(\ref{hF}) commute with respect to the Lie-Poisson bracket (\ref{hP}) on the dual of Lie algebra $so(2,1)^*$. 

At the symplectic leaf $||q||^2=1, \langle M,q\rangle=\nu$ this gives a new integrable two-centre problem on the hyperbolic plane with constant magnetic field of charge $\nu$.
The same is true for the natural quantum versions $\hat H$ and $\hat F.$
\end{theorem}

The formulae in the corresponding hyperbolic elliptic coordinates $u_1,u_2$, defined as the roots
$$\frac{q_1^2}{A-u}+\frac{q_2^2}{B-u}+\frac{q_3^2}{u}=0, \quad ||q||^2=-q_1^2-q_2^2+q_3^2=1,$$
are similar to the spherical case.

Note that on the symplectic leaves with $||q||^2=-1,$ determining a one-sheeted hyperboloid, we have a ``de Sitter" version of the problem using the same formulae with $B>0>A.$

\section{Concluding remarks}

% It would be interesting to find a geometric meaning of the potential of these new systems.

The geometry of the new systems is still to be properly studied, but the most important problem is to study the corresponding dynamics in the classical case and the spectrum in the quantum case. 

For the detailed analysis of the spherical version of the Euler two-centre problem \cite{K,KH} we refer to the work of Albouy and Stuchi \cite{Alb,AS}, Borisov and Mamaev \cite{BM,BM2} and Gonzalez Leon et al \cite{GGM}.

In our case this looks much more difficult because of the presence of the magnetic field, which usually creates a lot of problem for the separation of variables.

A famous example of this sort is the classical Clebsch system, describing special integrable cases of rigid body motion in an ideal fluid \cite{Clebsch}. It can be written as Euler equation on $e(3)^*$ with
 $$
 H=\frac{1}{2}|M|^2+\frac{1}{2}(Aq_1^2+Bq_2^2+Cq_3^2),
 $$
(see \cite{NS}) and can be interpreted as the harmonic oscillator on sphere with additional Dirac magnetic field \cite{V}. 
As we have already mentioned, Clebsch and our new systems are the only two electric extensions of Dirac magnetic monopole having an additional integral, which is quadratic in momenta \cite{VY}.
  
Recently there was a substantial progress in separation of variables for the Clebsch systems due to Magri and Skrypnyk \cite{MagriS, S}. It would be interesting to see if these new ideas can be applied in our case as well.

\section{Acknowledgements}

We are very grateful to Alexey Bolsinov for many helpful and stimulating discussions, and to Ivan Mamaev and Andrey Tsiganov for useful comments.

One of us (APV) would like to express special gratitude to Professor Valery Vasilievich Kozlov, who has just turned 70. He was my lecturer in theoretical mechanics at Moscow State University and his tutorials and papers, together with remarkable Arnold's book \cite{Arnold}, developed my taste for classical mechanics.

\end{document}